\UseRawInputEncoding
\documentclass[oneside, 12pt]{amsart}
\usepackage{amsmath}
\usepackage{amssymb}
\usepackage{dsfont}
\usepackage{color}
\usepackage[usenames,dvipsnames,x11names,svgnames]{xcolor}
\usepackage[colorlinks=true,linkcolor=NavyBlue,citecolor=DarkGreen, urlcolor=Blue]{hyperref}
\usepackage{amsthm}
\usepackage{amscd}
\usepackage{geometry}
\usepackage{mathrsfs}
\usepackage{mathtools}
\usepackage{esint}

\theoremstyle{definition}

\newcommand{\be}{\begin{equation*}}
\newcommand{\ee}{\end{equation*} }

\newcommand{\ben}{\begin{equation}}
\newcommand{\een}{\end{equation} }

\newcommand{\bs}{\begin{split}}
\newcommand{\es}{\end{split}}

\newcommand{\bmu}{\begin{multline*}}
\newcommand{\emu}{\end{multline*}}

\newcommand{\bmun}{\begin{multline}}
\newcommand{\emun}{\end{multline}}

\begin{document}



\title[]{An insight on the fractal power law flow: from a Hausdorff vector
calculus perspective}

\author{Xiao-Jun Yang}

\email{dyangxiaojun@163.com; xjyang@cumt.edu.cn}

\address{School of Mathematics, China University of Mining and Technology, Xuzhou 221116, China}

\begin{abstract}
In the article we suggest the Hausdorff vector calculus based on the Chen
Hausdorff calculus for the first time. The Gauss-Ostrogradsky-like,
Stokes-like, and Green-like theorems, and Green-like identities are obtained
in the framework of the Hausdorff vector calculus. The formula is proposed
as a mathematical tool to describe the real world problems for the fractal
power-law flow equations with the anomalous diffusion equation.
A conjecture for the fractal power-law flow equations analogous
to the Smales 15th Problem
(one of the Millennium Prize Problems for the Navier--Stokes equations) is also addressed.
\end{abstract}

\maketitle

\section{Introduction}\label{Sec:1}

The Hausdorff derivative involving the fractal geometry with the Hausdorff
measure, proposed by Chinese mathematician Wen Chen \cite{1,2}, has played an
important role in the treatment for the mathematical model for the anomalous
diffusion process.
The Hausdorff integral was suggested in 2018 by Chen and coauthors to
develope the three-dimensional diffusion model for fractal porous media \cite{3}.
The Hausdorff calculus was used to the mathematical models in the real world
problems for the fractal power law \cite{4}. As is known, the Chen Hausdorff
calculus is connected with the anomalous transport in porous media \cite{5} and
the fractional calculus \cite{6}.

The fractal vector calculus involving the fractal media was considered by
many researchers from the different~perspectives. The fractal vector
calculus involving the local fractional calculus was suggested in 2009 and
published in 2012 \cite{7}. The fractal vector calculus with the fractional
calculus was suggested \cite{8}. The vector calculus with the fractal metric
by using the fractional integrals was used to handle the fractional wave
problem \cite{9}. The fractal vector calculus was considered to model the fluid
flows in fractally permeable reservoirs by using the theory of the Hausdorff
derivative involving the Euclidean volume of the pore space scales \cite{10}.

The vector calculus with respect to monotone functions based on the Leibniz
derivative and Stieltjes-Riemann integral \cite{11,12} was proposed in 2020 by
author to suggest the PDEs arising in heat conduction \cite{13} and fluid flows
\cite{14}. By using the vector calculus with respect to monotone increasing
function, the vector power-law calculus was considered to deduce the PDEs in
the power-law fluid flow \cite{15}. The vector power-law calculus is called the
Hausdorff vector calculus if the power-laws has the same value. The main
target of the paper is to suggest the theory of the Hausdorff vector
calculus to suggest the real world problems for the fractal power-law flow
with the anomalous diffusion equation by using the Chen Hausdorff calculus.
The structure of the paper is
designed as follows. In Section \ref{Sec:2}, we introduce the theory of the Chen
Hausdorff calculus. In Section \ref{Sec:3}, we propose the theory of the Hausdorff
vector calculus. In Section \ref{Sec:4}, we present the system of the partial
differential equations arising in fractal power-law flow. In Section \ref{Sec:5}, we
suggest the anomalous diffusion equation arising in the real theory of the
turbulent fluid motion. Finally, we draw the conclusion in Section \ref{Sec:6}.

\section{The theory of the Chen Hausdorff calculus} \label{Sec:2}

Let ${\rm R}$ be the set of the real numbers.

Suppose that $\varpi \left( t \right)=t^\mu $ is the power-law function,
which derived from the Hausdorff measure \cite{1}, where $0<\mu \leq 1$ is the
fractal dimension, and $t\in \mathbb{R}$.

Let
\begin{equation}
\label{eq1}
\Phi \left( t \right)=\left( {\Phi \circ \varpi } \right)\left( t
\right)=\Phi \left( {\varpi \left( t \right)} \right),
\end{equation}
where $\Phi \left( \varpi \right)$ is the background function, which is
continuous.

We denote the sets of the composite functions by
\begin{equation}
\label{eq2}
\Im =\left\{ {\Phi \left( t \right):\Phi \left( t \right)=\Phi \left(
{\varpi \left( t \right)} \right),\varpi \left( t \right)=t^\mu } \right\}.
\end{equation}
Let $\Xi \in \Im $.

\begin{itemize}
\item \textbf{The Chen Hausdorff derivative}
\end{itemize}
The Chen Hausdorff derivative of the function $\Xi \left( t \right)$ is
defined as \cite{1}
\begin{equation}
\label{eq3}
{ }^CD_t^{\left({1} \right)} \Xi \left( t \right)=\frac{t^{1-\mu
}}{\mu }\frac{d\Xi \left( t \right)}{dt}.
\end{equation}
If (\ref{eq3}) holds, then $\Xi \left( t \right)$ is the Hausdorff differential
function for $t\in \mathbb{R}$, i.e.,
\[
\wp =\left\{ {g\left( t \right):g\left( t \right)={ }^CD_t^{\left( {1}
\right)} \Xi \left( t \right)} \right\}.
\]
The properties for the Chen Hausdorff derivative is given as follows:

(A1) The sum and difference rules for the Chen Hausdorff derivative:
\begin{equation}
\label{eq4}
{ }^CD_t^{\left( {1} \right)} \left( {\Xi _1 \left( t \right)\pm \Xi _2
\left( t \right)} \right)={ }^CD_t^{\left( {1} \right)} \Xi _1 \left( t
\right)\pm { }^CD_t^{\left( {1} \right)} \Xi _2 \left( t \right),
\end{equation}
where $\Xi _1 \in \wp $ and $\Xi _2 \in \wp $.

(A2) The constant multiple rule for the Chen Hausdorff derivative:
\begin{equation}
\label{eq5}
{ }^CD_t^{\left({1} \right)} \left( {\alpha \Xi \left( t \right)}
\right)=\alpha { }^CD_t^{\left({1} \right)} \Xi \left( t \right),
\end{equation}
where $\Xi \in \wp $ exists and $\alpha $ is a constant.

(A3) The product rule for the Chen Hausdorff derivative:
\begin{equation}
\label{eq6}
\begin{array}{l}
{ }^CD_t^{\left( {1} \right)} \left( {\Xi _1 \left( t \right)\cdot \Xi
_2 \left( t \right)} \right)
=\Xi _2 \left( t \right){ }^CD_t^{\left(
{1} \right)} \Xi _1 \left( t \right)+\Xi _1 \left( t \right){
}^CD_t^{\left( {1} \right)} \Xi _2 \left( t \right),
\end{array}
\end{equation}
where $\Xi _1 \in \wp $ and $\Xi _2 \in \wp $.

(A4) The quotient rule for the Chen Hausdorff derivative:
\begin{equation}
\label{eq7}
\begin{array}{l}
{ }^CD_t^{\left({1} \right)} \left( {\frac{\Xi _1 \left( t
\right)}{\Xi _2 \left( t \right)}} \right)
=\frac{\Xi _2 \left( t \right){
}^CD_t^{\left( {1} \right)} \Xi _1 \left( t \right)-\Xi _1 \left( t
\right){ }^CD_t^{\left( {1} \right)} \Xi _2 \left( t \right)}{\Xi _2
\left( t \right)\cdot \Xi _2 \left( t \right)},
\end{array}
\end{equation}
where $\Xi _1 \in \wp $ and $\Xi _2 \in \wp $ exist, and $\Xi _2 \left( t
\right)\ne 0$.

(A5) The chain rule for the Chen Hausdorff derivative:
\begin{equation}
\label{eq8}
{ }^CD_t^{\left({1} \right)} \Lambda \left( t \right)=\frac{d\Lambda
\left( \Xi \right)}{d\Xi }\cdot { }^CD_t^{\left( {1} \right)} \Xi
\left( t \right),
\end{equation}
where $\Lambda \left( t \right)=\Lambda \left( {\Xi \left( t \right)}
\right)=\left( {\Lambda \circ \Xi } \right)\left( t \right)$, $d\Lambda
\left( \Xi \right)/d\Xi $ and $\Xi \in \wp $.

The properties for the Chen Hausdorff derivative is presented as follows:
\begin{equation}
\label{eq9}
{ }^CD_t^{\left( {1} \right)} 1=0,
\end{equation}
\begin{equation}
\label{eq10}
{ }^CD_t^{\left( {1} \right)} \varpi \left( t \right)=1,
\end{equation}
\begin{equation}
\label{eq11}
{ }^CD_t^{\left( {1} \right)} \varpi ^n\left( t \right)=n\varpi
^{n-1}\left( t \right)=nt^{\mu \left( {n-1} \right)},
\end{equation}
\begin{equation}
\label{eq12}
{ }^CD_t^{\left( {1} \right)} e^{\beta t^\mu }=\beta e^{\beta t^\mu },
\end{equation}
\begin{equation}
\label{eq13}
{ }^CD_t^{\left({1} \right)} \ln \left( {t^\mu }
\right)=\frac{1}{t^\mu },
\end{equation}
\begin{equation}
\label{eq14}
{ }^CD_t^{\left( {1} \right)} s^{t^\mu }=\left( {\ln s} \right)s^{t^\mu
},
\end{equation}
\begin{equation}
\label{eq15}
{ }^CD_t^{\left( {1} \right)} \log _s \left( {t^\mu }
\right)=\frac{1}{t^\mu \ln s},
\end{equation}
\begin{equation}
\label{eq16}
{ }^CD_t^{\left( {1} \right)} e^{\Xi \left( t \right)}=e^{\Xi \left( t
\right)}{ }^CD_t^{\left( {1} \right)} \Xi \left( t \right),
\end{equation}
\begin{equation}
\label{eq17}
{ }^CD_t^{\left( {1} \right)} \log _s \Xi \left( t \right)=\frac{1}{\ln
s}\frac{{ }^CD_t^{\left( {1} \right)} \Xi \left( t \right)}{\Xi \left(
t \right)},
\end{equation}
\begin{equation}
\label{eq18}
{ }^CD_t^{\left( {1} \right)} \ln \Xi \left( t \right)=\frac{{
}^CD_t^{\left( {1} \right)} \Xi \left( t \right)}{\Xi \left( t
\right)},
\end{equation}
\begin{equation}
\label{eq19}
{ }^CD_t^{\left( {1} \right)} s^{\Xi \left( t \right)}=\left[ {\left(
{\ln s} \right)s^{\Xi \left( t \right)}} \right]\cdot {\begin{array}{*{20}c}
 {{ }^CD_t^{\left( {1} \right)} \Xi \left( t \right)} \hfill \\
\end{array} },
\end{equation}
where $\beta $ is the constant and
\begin{equation}
\label{eq20}
e^{\beta t^\mu }=\sum\limits_{n=0}^\infty {\left( {\beta ^nt^{n\mu }}
\right)/n!} ,
\end{equation}
is the Kohlrausch-Williams-Watts function \cite{12}.

\begin{itemize}
\item \textbf{The Chen Hausdorff integral}
\end{itemize}
Let $\xi \in \Im $.

The Chen Hausdorff integral of the function $\xi \left( t \right)$ in the
interval $\left[ {a,b} \right]$ is defined as \cite{3}
\begin{equation}
\label{eq21}
{ }_a^C I_b^{\left( {1} \right)} \xi \left( t \right)=\mu
\int\limits_a^b {\xi \left( t \right)t^{\mu -1}dt} .
\end{equation}
If (\ref{eq3}) holds, then $\Xi \left( t \right)$ is the Hausdorff integral in the
interval $\left[ {a,b} \right]$, i.e.,
\[
\aleph =\left\{ {h\left( t \right):h\left( t \right)=\mu \int\limits_a^b
{\xi \left( t \right)t^{\mu -1}dt} } \right\}.
\]
The properties for the Chen Hausdorff integral is presented as follows:

(B1) The sum and difference rules for the Chen Hausdorff integral:
\begin{equation}
\label{eq22}
{ }_a^C I_b^{\left( {1} \right)} \left( {\xi _1 \left( t \right)\pm \xi
_2 \left( t \right)} \right)={ }_a^C I_b^{\left( {1} \right)} \xi _1
\left( t \right)\pm { }_a^C I_b^{\left({1} \right)} \xi _2 \left( t
\right),
\end{equation}
where $\xi _1 \in \aleph $ and $\xi _2 \in \aleph $.

(B2) The first fundamental theorem of the Chen Hausdorff integral:
\begin{equation}
\label{eq23}
\Xi \left( t \right)-\Xi \left( a \right)={ }_a^C I_t^{\left( {1}
\right)} \left( {{\begin{array}{*{20}c}
 {{ }^CD_t^{\left( {1} \right)} } \hfill \\
\end{array} }\Xi \left( t \right)} \right)
\end{equation}
(B3) The mean value theorem for the topology integral:
\begin{equation}
\label{eq24}
{ }_a^C I_t^{\left( {1} \right)} \xi \left( t \right)=\xi \left( l
\right)\left( {\varpi \left( t \right)-\varpi \left( a \right)} \right)
\end{equation}
where $a<t\le l<b$, and $\xi \in \aleph $.

(B4) The second fundamental theorem of the Chen Hausdorff integral:
\begin{equation}
\label{eq25}
\xi \left( t \right)={\begin{array}{*{20}c}
 {{ }^CD_t^{\left( {1} \right)} } \hfill \\
\end{array} }\left( {{ }_a^C I_t^{\left( {1} \right)} \xi \left( t
\right)} \right)
\end{equation}
where $\xi \in \aleph $.

(B5) The net change theorem for the Chen Hausdorff integral:
\begin{equation}
\label{eq26}
\Xi \left( b \right)-\Xi \left( a \right)={ }_a^T I_b^{\left( {1}
\right)} \left( {{\begin{array}{*{20}c}
 {{ }^TD_t^{\left( {1} \right)} } \hfill \\
\end{array} }\Xi \left( t \right)} \right)
\end{equation}
where $\Xi \in \wp $.

(B6) The integration by parts for the Chen Hausdorff integral:
\begin{equation}
\label{eq27}
\begin{array}{l}
{ }_a^C I_b^{\left( {1} \right)} \left( {\xi _2 \left( t \right){
}^CD_t^{\left( {1} \right)} \xi _1 \left( t \right)}
\right)
=\left[ {\xi _1 \left( t \right)\cdot \xi _2 \left( t \right)}
\right]_a^b -{ }_a^C I_b^{\left( {1} \right)} \left( {\xi _1 \left( t
\right){ }^CD_t^{\left( {1} \right)} \xi _2 \left( t \right)} \right),
\end{array}
\end{equation}
where $\left[ {\xi _1 \left( t \right)\cdot \xi _2 \left( t \right)}
\right]_a^b =\xi _1 \left( b \right)\cdot \xi _2 \left( b \right)-\xi _1
\left( a \right)\cdot \xi _2 \left( a \right)$, $\xi _1 \in \aleph $ and
$\xi _2 \in \aleph $.

The indefinite Chen Hausdorff integral of the function $\xi \left( t
\right)$ is defined as
\[
{ }^CI_t^{\left( {1} \right)} \xi \left( t \right)=\mu \int {\xi \left(
t \right)t^{\mu -1}dt} =\Xi \left( t \right){+}\chi ,
\]
which implies that
\[
{ }^CD_t^{\left({1} \right)} \left( {{ }^CI_t^{\left( {1}
\right)} \xi \left( t \right)} \right)={ }^CD_t^{\left( {1} \right)}
\Xi \left( t \right)=\xi \left( t \right)
\]
and
\[
\begin{array}{l}
\Xi \left( t \right)-\Xi \left( a \right)
={ }^CI_t^{\left( {1} \right)}
\xi \left( t \right)-{ }^CI_a^{\left( {1} \right)} \xi \left( t
\right)
={ }_a^C I_t^{\left( {1} \right)} \left( {{\begin{array}{*{20}c}
 {{ }^CD_t^{\left( {1} \right)} } \hfill \\
\end{array} }\Xi \left( t \right)} \right),
\end{array}
\]
where $\chi $ is the constant.

The sum and difference rules for the indefinite Chen Hausdorff integral
implies that
\[
{ }^CI_t^{\left( {1} \right)} \left( {\xi _1 \left( t \right)\pm \xi _2
\left( t \right)} \right)={ }^CI_t^{\left( {1} \right)} \xi _1 \left( t
\right)\pm { }^CI_t^{\left( {1} \right)} \xi _2 \left( t \right),
\]
where $\xi _1 \in \aleph $ and $\xi _2 \in \aleph $.

The properties for the Chen Hausdorff integral is presented as follows:
\begin{equation}
\label{eq28}
{ }^CI_t^{\left( {1} \right)} 1=\varpi \left( t \right){+}\chi
=t^\mu {+}\chi ,
\end{equation}
\begin{equation}
\label{eq29}
{ }^CI_t^{\left({1} \right)} \left( {n\varpi ^{n-1}\left( t \right)}
\right)=\varpi ^n\left( t \right)+\chi ,
\end{equation}
\begin{equation}
\label{eq30}
{ }^CI_t^{\left( {1} \right)} \left( {\frac{1}{\ln s}\cdot \frac{{
}^CD_t^{\left( {1} \right)} \Xi \left( t \right)}{\Xi \left( t
\right)}} \right)=\log _s \Xi \left( t \right)+\chi ,
\end{equation}
\begin{equation}
\label{eq31}
{ }^CI_t^{\left( {1} \right)} \frac{1}{\varpi \left( t \right)}=\ln
\varpi \left( t \right)+\chi ,
\end{equation}
\begin{equation}
\label{eq32}
{ }^CI_t^{\left( {1} \right)} \left( {\frac{1}{\ln s}\cdot
\frac{1}{\varpi \left( t \right)}} \right)=\log _s \varpi \left( t
\right)+\chi ,
\end{equation}
\begin{equation}
\label{eq33}
{ }^CI_t^{\left( {1} \right)} \left( {\left( {\ln s} \right)s^{t^\mu }}
\right)=s^{t^\mu }+\chi ,
\end{equation}
\begin{equation}
\label{eq34}
{ }^CI_t^{\left( {1} \right)} \left( {e^{\Xi \left( t \right)}{
}^CD_t^{\left( {1} \right)} \Xi \left( t \right)} \right)=e^{\Xi \left(
t \right)}+\chi ,
\end{equation}
\begin{equation}
\label{eq35}
{ }^CI_t^{\left( {1} \right)} \left( {\frac{\Xi \left( t
\right)}{\left| {\Xi \left( t \right)} \right|}{ }^CD_t^{\left( {1}
\right)} \Xi \left( t \right)} \right)=\left| {\Xi \left( t \right)}
\right|+\chi ,
\end{equation}
\begin{equation}
\label{eq36}
{ }^CI_t^{\left( {1} \right)} \left( {\frac{{ }^CD_t^{\left( {1}
\right)} \Xi \left( t \right)}{\Xi \left( t \right)}} \right)=\ln \Xi \left(
t \right)+\chi ,
\end{equation}
\begin{equation}
\label{eq37}
{ }^CI_t^{\left( {1} \right)} \left( {e^{\beta t^\mu }} \right)=\beta
e^{\beta t^\mu }+\chi ,
\end{equation}
\begin{equation}
\label{eq38}
{ }^CI_t^{\left( {1} \right)} \left[ {\left( {\ln s} \right)s^{\Xi
\left( t \right)}\cdot { }^CD_t^{\left( {1} \right)} \Xi \left( t
\right)} \right]=s^{\Xi \left( t \right)}+\chi ,
\end{equation}
where $\chi $ is the constant.

\begin{itemize}
\item \textbf{The Chen Hausdorff partial derivatives}
\end{itemize}
We now consider the coordinate system, expressed by
\begin{equation}
\label{eq39}
{\rm {\bf i}}x^\mu +{\rm {\bf j}}y^\mu +{\rm {\bf k}}z^\mu =\left( {x^\mu
,y^\mu ,z^\mu } \right),
\end{equation}
where ${\rm {\bf i}}$, ${\rm {\bf j}}$ and ${\rm {\bf k}}$ denote the unit
vectors in the Cartesian coordinate system.

We now define the function by
\begin{equation}
\label{eq40}
\psi =\Pi \left( {x^\mu ,y^\mu ,z^\mu } \right).
\end{equation}
The Chen partial derivatives of the function $\psi =\Pi \left( {x^\mu ,y^\mu
,z^\mu } \right)$ are defined as
\begin{equation}
\label{eq41}
{ }^C\partial _x^{\left( 1 \right)} \psi =\left( {\frac{x^{1-\mu }}{\mu
}\frac{\partial }{\partial x}} \right)\psi ,
\end{equation}
\begin{equation}
\label{eq42}
{ }^C\partial _y^{\left( 1 \right)} \psi =\left( {\frac{x^{1-\mu }}{\mu
}\frac{\partial }{\partial x}} \right)\psi ,
\end{equation}
\begin{equation}
\label{eq43}
{ }^C\partial _z^{\left( 1 \right)} \psi =\left( {\frac{x^{1-\mu }}{\mu
}\frac{\partial }{\partial x}} \right)\psi .
\end{equation}
The total differential of the function $\psi =\Pi \left( {x^\mu ,y^\mu
,z^\mu } \right)$ is defined as
\begin{equation}
\label{eq44}
\begin{array}{l}
d\psi =\mu {\left( {x^{\mu -1}{ }^C\partial _x^{\left( 1 \right)}
\psi } \right)dx}
       +\mu {\left( {y^{\mu -1}{ }^C\partial _y^{\left( 1 \right)} \psi
} \right)dy}\\
       +\mu {\left( {z^{\mu -1}{ }^C\partial _z^{\left( 1 \right)} \psi }
\right)dz}.
\end{array}
\end{equation}
Thus, the Chen Hausdorff derivative with respect to the time $t$ reads
\begin{equation}
\label{eq45}
\begin{array}{l}
\frac{d\psi }{dt}
=\mu \left( {x^{\mu -1}{ }^C\partial _x^{\left( 1 \right)} \psi } \right)\frac{dx}{dt}+
 \mu \left( {y^{\mu -1}{ }^C\partial _y^{\left( 1 \right)} \psi } \right)\frac{dy}{dt}\\
+\mu \left( {z^{\mu -1}{ }^C\partial _z^{\left( 1 \right)} \psi } \right)\frac{dz}{dt}.
\end{array}
\end{equation}
\begin{itemize}
\item \textbf{The Chen gradient}
\end{itemize}
The Chen gradient in the Cartesian coordinate system is defined as \cite{1,4}
\begin{equation}
\label{eq46}
\begin{array}{l}
\nabla ^\mu
=\mu {\rm {\bf i}}\left( {x^{\mu -1}} \right){}^C\partial _x^{\left( 1 \right)}
+{\rm {\bf j}}\left( {y^{\mu -1}} \right){}^C\partial _y^{\left( 1 \right)}
+{\rm {\bf k}}\left( {z^{\mu -1}} \right){}^C\partial _z^{\left( 1 \right)} .
\end{array}
\end{equation}
From (\ref{eq44}) we arrive at
\begin{equation}
\label{eq47}
d\psi =\nabla ^\mu \psi \cdot d{\rm {\bf r}}=\nabla ^\mu \psi \cdot {\rm
{\bf n}}dr,
\end{equation}
in which
\begin{equation}
\label{eq48}
d{\rm {\bf r}}={\rm {\bf n}}dr={\rm {\bf i}}dx+{\rm {\bf j}}dy+{\rm {\bf
k}}dz,
\end{equation}
where ${\rm {\bf n}}$ is the unit normal, and $dr$ is a distance measured
along the normal ${\rm {\bf n}}$.

\begin{itemize}
\item \textbf{The Hausdorff directional derivative }
\end{itemize}
The Hausdorff directional derivative of the function $\psi =\Pi \left(
{x^\mu ,y^\mu ,z^\mu } \right)$, denoted by $\nabla _n^\mu \psi $, is
defined as
\begin{equation}
\label{eq49}
\frac{d\psi }{dr}=\nabla ^\mu \psi \cdot {\rm {\bf n}}=\partial _n^\mu \psi
,
\end{equation}
where $d\psi /dr$ is the rates of change of $\psi $ along the normal ${\rm
{\bf n}}$, respectively.

\begin{itemize}
\item \textbf{The Laplace-Chen operator }
\end{itemize}
The Laplace-Chen operator, denoted as $\nabla ^\mu \cdot \nabla ^\mu =\nabla
^{2\mu }$, of the scalar field $\psi $ is defined as \cite{1,4}
\begin{equation}
\label{eq50}
\begin{array}{l}
 \nabla ^{2\mu }\psi \\
 =\mu ^2 {\left( {x^{\mu -1}{ }^C\partial_x^{\left( 1 \right)} } \right)^2}\psi
 +\mu ^2 {\left( {y^{\mu -1}{ }^C\partial _y^{\left(1 \right)} } \right)^2}\psi \\
 +\mu ^2 {\left( {z^{\mu -1}{ }^C\partial _z^{\left( 1 \right)}} \right)^2}\psi \\
 =\mu ^2 {{x^{2\mu -2}{ }^C\partial _x^{\left( 2 \right)} } }\psi
 +\mu ^2 {{y^{2\mu -2}{ }^C\partial _y^{\left( 2 \right)} } }\psi\\
 +\mu ^2 {{z^{2\mu -2}{ }^C\partial _z^{\left( 2 \right)} } }\psi. \\
 \end{array}
\end{equation}
The properties for the Laplace-Chen operator are presented as
\begin{equation}
\label{eq51}
\left( {\nabla ^\mu \cdot \nabla ^\mu } \right)\psi =\nabla ^{2\mu }\psi ,
\end{equation}
\begin{equation}
\label{eq52}
\nabla ^\mu \left( {\psi \Theta } \right)=\psi \nabla ^\mu \Theta +\Theta
\nabla ^\mu \psi ,
\end{equation}
\begin{equation}
\label{eq53}
\nabla ^\mu \cdot \left( {\Theta \nabla ^\mu \psi } \right)=\Theta \nabla
^{2\mu }\psi +\nabla ^\mu \psi \cdot \nabla ^\mu \Theta ,
\end{equation}
where $\psi $ and $\Theta $ are the fractal scalar fields.

Here, (\ref{eq51}) is discovered by Chen \cite{1}.

\section{The theory of the Hausdorff vector calculus}\label{Sec:3}

The element of the vector line
\begin{equation}
\label{eq54}
\ell =\ell \left( {x,y,z} \right)=\widetilde{\ell }\left( {x^\mu ,y^\mu
,z^\mu } \right)
\end{equation}
is expressed in the form
\begin{equation}
\label{eq55}
\begin{array}{l}
d{\rm {\bf l}}
={\rm {\bf m}}d\ell \\
=\mu \left( {{\rm {\bf i}}x^{\mu
-1}dx+{\rm {\bf j}}y^{\mu -1}dy+{\rm {\bf k}}z^{\mu -1}dz} \right)
\end{array}
\end{equation}
and
\begin{equation}
\label{eq56}
\begin{array}{l}
d\ell =\left| {d{\rm {\bf l}}} \right|\\
=\mu \sqrt {x^{2\mu -2}\left( {dx}
\right)^2+y^{2\mu -2}\left( {dy} \right)^2+z^{2\mu -2}\left( {dz} \right)^2}
,
\end{array}
\end{equation}
where ${\rm {\bf m}}$ is the vector with $\left| {\rm {\bf m}} \right|=1$.

The Hausdorff arc length is represented in the form:
\begin{equation}
\label{eq57}
\begin{array}{l}
\ell= \\
\int\limits_a^b {\sqrt {x^{2\mu -2}\left( {\frac{dx}{dt}}
\right)^2+y^{2\mu -2}\left( {\frac{dy}{dt}} \right)^2+z^{2\mu -2}\left(
{\frac{dz}{dt}} \right)^2} } dt.
\end{array}
\end{equation}
\begin{itemize}
\item \textbf{The line Hausdorff integral of the fractal vector field}
\end{itemize}
The line Hausdorff integral of the fractal vector field
\[
{\rm {\bf {\rm T}}}={\rm {\bf {\rm T}}}\left( {x,y,z}
\right)=\widetilde{{\rm {\bf {\rm T}}}}\left( {x^\mu ,y^\mu ,z^\mu }
\right)
\]
along the curve
\[
L=L\left( {x,y,z} \right)=\widetilde{L}\left( {x^\mu ,y^\mu ,z^\mu }
\right),
\]
denoted by ${\rm L}$, is defined as
\begin{equation}
\label{eq58}
{\rm L}=\int\limits_L {{\rm {\bf {\rm T}}}\left( {x,y,z} \right)\cdot d{\rm
{\bf l}}} ,
\end{equation}
where the element of the vector line is
\begin{equation}
\label{eq59}
d{\rm {\bf l}}=\mu \left( {{\rm {\bf i}}x^{\mu -1}dx+{\rm {\bf j}}y^{\mu
-1}dy+{\rm {\bf k}}z^{\mu -1}dz} \right).
\end{equation}
By using (\ref{eq62}), we show that
\begin{equation}
\label{eq60}
\int\limits_L {{\rm {\bf {\rm T}}}\cdot d{\rm {\bf l}}} =\int\limits_L {{\rm
{\bf {\rm T}}}\left( {x,y,z} \right)\cdot d{\rm {\bf l}}}
=\int\limits_{L\left( t \right)} {{\rm {\bf {\rm T}}}\cdot \frac{d{\rm {\bf
l}}}{dt}dt} ,
\end{equation}
where
\[
d{\rm {\bf l}}/dt=\mu \left( {{\rm {\bf i}}x^{\mu -1}dx/dt+{\rm {\bf
j}}y^{\mu -1}dy/dt+{\rm {\bf k}}z^{\mu -1}dz/dt} \right).
\]
Thus, from (\ref{eq58}), (\ref{eq58}) can be presented as follows:
\begin{equation}
\label{eq61}
\int\limits_L {{\rm {\bf {\rm T}}}\cdot d{\rm {\bf l}}} =\mu \int\limits_L
{T_x x^{\mu -1}dx+T_y y^{\mu -1}dy+T_z z^{\mu -1}dz} .
\end{equation}
The vector field ${\rm {\bf {\rm T}}}={\rm {\bf {\rm T}}}\left( {x,y,z}
\right)$ in
\[
L=L\left( {x,y,z} \right)=\widetilde{L}\left( {x^\mu ,y^\mu ,z^\mu }
\right)
\]
is said to be conservative if
\begin{equation}
\label{eq62}
\oint\limits_L {{\rm {\bf {\rm T}}}\cdot d{\rm {\bf l}}} =0.
\end{equation}
\begin{itemize}
\item \textbf{The double Hausdorff integral of the fractal scalar field}
\end{itemize}
The double Hausdorff integral of the fractal scalar field
\[
{\rm M}={\rm M}\left( {x,y} \right)=\widetilde{{\rm M}}\left( {x^\mu ,y^\mu
} \right)
\]
on the region $S\left( {x,y} \right)=\widetilde{S}\left( {x^\mu ,y^\mu }
\right)$, denoted by $A\left( {\rm M} \right)$, is defined as
\begin{equation}
\label{eq63}
A\left( {\rm M} \right)=\mathop{\iint}\limits_S {{\rm M}\left( {x,y}
\right)dS} ,
\end{equation}
where
\begin{equation}
\label{eq64}
dS=\mu ^2x^{\mu -1}y^{\mu -1}dxdy.
\end{equation}
With the aid of from (\ref{eq63}) and (\ref{eq64}) we may see that
\begin{equation}
\label{eq65}
\begin{array}{l}
\mathop{\iint}\limits_S {{\rm M}\left( {x,y} \right)dS}\\
 =\mu
^2\int\limits_c^d {\left[ {\int\limits_a^b {{\rm M}\left( {x,y}
\right)x^{\mu -1}dx} } \right]y^{\mu -1}dy} \\
=\mu ^2\int\limits_a^b {\left[ {\int\limits_c^d {{\rm M}\left(
{x,y} \right)y^{\mu -1}dy} } \right]x^{\mu -1}dx} , \\
 \end{array}
\end{equation}
where $x\in \left[ {a,b} \right]$ and $y\in \left[ {c,d} \right]$.

\begin{itemize}
\item \textbf{The volume Hausdorff integral of the fractal scalar field}
\end{itemize}
The volume Hausdorff integral of the fractal scalar field
\[
N=N\left( {x,y,z} \right)=\widetilde{N}\left( {x^\mu ,y^\mu ,z^\mu }
\right)
\]
is defined as
\begin{equation}
\label{eq66}
V\left( N \right)=\mathop{\iiint}\limits_{\kern-5.5pt
\Omega } {N\left( {x,y,z} \right)dV} ,
\end{equation}
where
\[
dV=\mu ^3x^{\mu -1}y^{\mu -1}z^{\mu -1}dxdydz
\]
and
\[
\Omega =\Omega \left( {x,y,z} \right)=\widetilde{\Omega }\left( {x^\mu
,y^\mu ,z^\mu } \right).
\]
By using (\ref{eq66}), we may show
\begin{equation}
\label{eq67}
\begin{array}{l}
\mathop{\iiint} \limits_{\kern-5.5pt \Omega } {N\left(
{x,y,z} \right)dV} \\
=\mu ^3\int\limits_\alpha ^\beta {z^{\mu
-1}dz\int\limits_c^d {y^{\mu -1}dy\int\limits_a^b {N\left( {x,y,z}
\right)x^{\mu -1}dx} } } \\
=\mu ^3\int\limits_a^b {x^{\mu -1}dx\int\limits_\alpha ^\beta
{z^{\mu -1}dz\int\limits_c^d {N\left( {x,y,z} \right)y^{\mu -1}dy} } } \\
=\mu ^3\int\limits_c^d {y^{\mu -1}dy\int\limits_a^b {x^{\mu
-1}dx\int\limits_\alpha ^\beta {Nz^{\mu -1}dz} } }, \\
 \end{array}
\end{equation}
where $x\in \left[ {a,b} \right]$, $y\in \left[ {c,d} \right]$ and $z\in
\left[ {\alpha ,\beta } \right]$.

\begin{itemize}
\item \textbf{The surface Hausdorff integral of the fractal vector field}
\end{itemize}
The surface Hausdorff integral of the fractal vector field
\[
\begin{array}{l}
{\rm {\bf W}}={\rm {\bf W}}\left( {x,y,z} \right)=\widetilde{{\rm {\bf
W}}}\left( {x^\mu ,y^\mu ,z^\mu } \right)
 \end{array}
\]
is defined as
\begin{equation}
\label{eq68}
\begin{array}{l}
\mathop{\iint} \limits_{\rm {\bf S}} {{\rm {\bf W}}\left( {x,y,z}
\right)\cdot d{\rm {\bf S}}} =\mathop{\iint} \limits_{\rm {\bf S}} {{\rm {\bf
W}}\left( {x,y,z} \right)\cdot {\rm {\bf n}}dS} ,
 \end{array}
\end{equation}
where ${\rm {\bf n}}=d{\rm {\bf S}}/dS$ is the unit normal vector to
the surface
\[
\begin{array}{l}
{\rm {\bf S}}={\rm {\bf S}}\left( {x,y,z} \right)=\widetilde{{\rm {\bf
s}}}\left( {x^\mu ,y^\mu ,z^\mu } \right).
 \end{array}
\]
Suppose that
\[
\begin{array}{l}
{\rm {\bf n}}={\rm {\bf dS}}/\left| {{\rm {\bf dS}}} \right|={\rm {\bf
dS}}/dS,
 \end{array}
\]

\[
\begin{array}{l}
dS=\left| {{\rm {\bf dS}}} \right|,
 \end{array}
\]
and
\begin{equation}
\label{eq69}
\begin{array}{l}
{\rm {\bf dS}}=\mu ^2 {\rm {\bf i}}y^{\mu -1}z^{\mu -1}dydz\\
+\mu ^2{\rm {\bf j}}x^{\mu -1}z^{\mu -1}dxdz
+\mu ^2{\rm {\bf k}}x^{\mu -1}y^{\mu -1}dxdy.
\end{array}
\end{equation}
Then, by (\ref{eq69}), Eq.(\ref{eq68}) can be represented in the form:
\begin{equation}
\label{eq70}
\begin{array}{l}
\int\!\!\!\int\limits_{\rm {\bf S}} {{\rm {\bf W}}\left( {x,y,z}
\right)\cdot d{\rm {\bf S}}} \\
=\mu ^2 {\int\!\!\!\int\limits_S {W_x
y^{\mu -1}z^{\mu -1}dydz} + \mu ^2 W_y x^{\mu -1}z^{\mu -1}dxdz} \\
+\mu ^2 {W_z x^{\mu -1}y^{\mu
-1}dxdy}
\end{array}
\end{equation}
where
\[
\begin{array}{l}
{\rm {\bf W}}={\rm {\bf W}}\left( {x,y,z} \right)\\
=\widetilde{{\rm {\bf
w}}}\left( {x^\mu ,y^\mu ,z^\mu } \right)\\
={\rm {\bf i}}W_x +{\rm {\bf j}}W_y
+{\rm {\bf k}}W_z .
\end{array}
\]
The flux of the fractal vector field ${\rm {\bf W}}={\rm {\bf W}}\left(
{x,y,z} \right)$ across the surface $d{\rm {\bf S}}$, denoted by ${\rm {\bf
Q}}$, is defined as
\begin{equation}
\label{eq71}
{\rm {\bf Q}}=\oiint \limits_{\rm
{\bf S}} {{\rm {\bf W}}\left( {x,y,z} \right)\cdot d{\rm {\bf S}}} .
\end{equation}
Let ${\rm {\bf Q}}=0$. Then we have from (\ref{eq71}) that
\begin{equation}
\label{eq72}
\oiint \limits_{\rm {\bf S}} {{\rm
{\bf W}}\left( {x,y,z} \right)\cdot d{\rm {\bf S}}} =0.
\end{equation}
\begin{itemize}
\item \textbf{The Hausdorff divergence of the fractal vector field}
\end{itemize}
The Hausdorff divergence of the fractal vector field ${\rm {\bf \psi }}$ is
defined as
\begin{equation}
\label{eq73}
\nabla ^\mu \cdot {\rm {\bf W}}=\mathop {\lim }\limits_{\Delta V_m \to 0}
\frac{1}{\Delta V_m } \oiint \limits_{\Delta {\rm {\bf S}}_m } {{\rm {\bf W}}\cdot d{\rm {\bf
S}}} ,
\end{equation}
where the volume $V$ is divided into a large number of small subvolumes
$\Delta V_m $ with surfaces $\Delta {\rm {\bf S}}_m $, ${\rm {\bf W}}$ is a Hausdorff
differentiable vector field, and $d{\rm {\bf S}}$ is an element of the
surface ${\rm {\bf S}}$ bounding the solid $\Omega $.

In the coordinate system (\ref{eq39}), (\ref{eq73}) can be represented as \cite{1}
\begin{equation}
\label{eq74}
\begin{array}{l}
\nabla ^\mu \cdot {\rm {\bf W}}\\
=\mu \left( {x^{\mu -1}{ }^C\partial
_x^{\left( 1 \right)} W_x +y^{\mu -1}{ }^C\partial _y^{\left( 1 \right)} W_y
+z^{\mu -1}{ }^C\partial _z^{\left( 1 \right)} W_z } \right),
\end{array}
\end{equation}
where
\[
{\rm {\bf W}}={\rm {\bf W}}\left( {x,y,z} \right)=\widetilde{{\rm {\bf
W}}}\left( {x^\mu ,y^\mu ,z^\mu } \right)={\rm {\bf i}}W_x +{\rm {\bf j}}W_y
+{\rm {\bf k}}W_z .
\]
\begin{itemize}
\item \textbf{The Hausdorff curl of the fractal vector field}
\end{itemize}
The Hausdorff curl of the fractal vector field ${\rm {\bf W}}$ is defined as
\begin{equation}
\label{eq75}
\nabla ^\mu \times {\rm {\bf W}}=\mathop {\lim }\limits_{\Delta V_m \to 0}
\frac{1}{\Delta V_m } \oiint \limits_{\Delta {\rm {\bf S}}_m } {{\rm {\bf W}}\times d{\rm {\bf
S}}} ,
\end{equation}
where the volume $V$ is divided into a large number of small subvolumes
$\Delta V_m $ with surfaces $\Delta {\rm {\bf S}}_m $, ${\rm {\bf W}}$ is a Hausdorff
differentiable vector field, and $d{\rm {\bf S}}$ is an element of the
surface ${\rm {\bf S}}$ bounding the solid $\Omega $.

There is an alternative definition of (\ref{eq75}) as follows:

The Hausdorff curl of the fractal vector field ${\rm {\bf W}}$ is defined as
\begin{equation}
\label{eq76}
\left( {\nabla ^\mu \times {\rm {\bf W}}} \right)\cdot {\rm {\bf n}}=\mathop
{\lim }\limits_{\Delta S_m \to 0} \frac{1}{\Delta S_m
}\oint\limits_{\Delta L_m } {{\rm {\bf W}}\cdot d{\rm
{\bf l}}} ,
\end{equation}
where ${\rm {\bf W}}$ is a Hausdorff differentiable vector field, $d{\rm {\bf
l}}$ is the element of the vector line, $\Delta S_m $ is a small surface
element perpendicular to ${\rm {\bf n}}$, $\Delta L_m $ is the closed curve
of the boundary of $\Delta S_m $, and ${\rm {\bf n}}$ are oriented in a
postive sense.

Similarly, in the coordinate system (\ref{eq39}), Eqs.(\ref{eq75}) and (\ref{eq76}) can be rewritten
as \cite{1}
\begin{equation}
\label{eq77}
\begin{array}{l}
\nabla ^\mu \times {\rm {\bf W}}
=\left( {{\begin{array}{*{20}c}
 {\rm {\bf i}} \hfill & {\rm {\bf j}} \hfill & {\rm {\bf k}} \hfill \\
 {\mu x^{\mu -1}{ }^C\partial _x^{\left( 1 \right)} } \hfill & {\mu x^{\mu
-1}{ }^C\partial _y^{\left( 1 \right)} } \hfill & {\mu x^{\mu -1}{
}^C\partial _z^{\left( 1 \right)} } \hfill \\
 {W_x } \hfill & {W_y } \hfill & {W_z } \hfill \\
\end{array} }} \right),
\end{array}
\end{equation}
where
\[
\begin{array}{l}
{\rm {\bf W}}
={\rm {\bf W}}\left( {x,y,z} \right)\\
=\widetilde{{\rm {\bf
W}}}\left( {x^\mu ,y^\mu ,z^\mu } \right)={\rm {\bf i}}W_x +{\rm {\bf j}}W_y
+{\rm {\bf k}}W_z .
\end{array}
\]
\begin{itemize}
\item \textbf{The Gauss-Ostrogradsky-like theorem for the fractal vector field}
\end{itemize}
By using (\ref{eq73}), we obtain the Gauss-Ostrogradsky-like theorem for the fractal
vector field ${\rm {\bf W}}$.

Now, we show that
\begin{equation}
\label{eq78}
\begin{array}{l}
\oiint \limits_{\bf{S}}
 {{\bf{W}}{\bf{n}}dS}  = \oiint \limits_{\bf{S}}
 {{\bf{W}}d{\bf{S}}}  \\
 {\rm{                }} \\
  = \mu ^2 \oiint \limits_{\bf{S}}
 {\left( {{\bf{i}}W_x  + {\bf{j}}W_y  + {\bf{k}}W_z } \right)}  \\
 \cdot
 \left( {{\bf{i}}y^{\mu  - 1} z^{\mu  - 1} dydz + {\bf{j}}x^{\mu  - 1} z^{\mu  - 1} dxdz + {\bf{k}}x^{\mu  - 1} y^{\mu  - 1} dxdy} \right) \\
 {\rm{                }} \\
  = \mu ^2 \oiint \limits_{\bf{S}}
 {\left[ {W_x y^{\mu  - 1} z^{\mu  - 1} dydz + W_y x^{\mu  - 1} z^{\mu  - 1} dxdz} \right]}  \\
 + \mu ^2 \oiint \limits_{\bf{S}}
{W_z x^{\mu  - 1} y^{\mu  - 1} dxdy} . \\
 \end{array}
\end{equation}

From (\ref{eq73}) the Gauss-like theorem for the fractal vector field ${\rm {\bf
W}}$states that
\begin{equation}
\label{eq79}
\mathop{\iiint} \limits_{\kern-5.5pt \Omega } {\nabla ^\mu
\cdot {\rm {\bf W}}dV} =\oiint \limits_{\rm {\bf S}} {{\rm {\bf W}}\cdot {\rm {\bf n}}dS} ,
\end{equation}
where ${\rm {\bf W}}$ is a Hausdorff differentiable vector field,
$dV$ denotes an element of volume $\Omega $, ${\rm {\bf n}}$ is the unit
outward normal to ${\rm {\bf S}}$, and $dS$ is an element of the surface
area of the surface ${\rm {\bf S}}$ bounding the solid $\Omega $.

With (\ref{eq78}) we have
\[
d{\rm {\bf S}}={\rm {\bf n}}dS,
\]
we obtain an alternative form of (\ref{eq79}) as follows:
\begin{equation}
\label{eq80}
\mathop{\iiint} \limits_{\kern-5.5pt \Omega } {\nabla ^\mu
\cdot {\rm {\bf W}}dV} =\oiint \limits_{\rm {\bf S}} {{\rm {\bf W}}\cdot d{\rm {\bf S}}} .
\end{equation}
It is easy to see that (\ref{eq80}) becomes the Gauss-Ostrogradsky theorem \cite{16} due
to Gauss \cite{17} and Ostrogradsky \cite{18} if $\mu =1$.

\begin{itemize}
\item \textbf{The Stokes-like theorem for the fractal vector field}
\end{itemize}
By using (\ref{eq76}), we present the Stokes-like theorem for the fractal vector
field.

We now consider that
\begin{equation}
\label{eq81}
\begin{array}{l}
 \oint\limits_L {{\bf{W}}d{\bf{l}}} \\
  = \mu \oint\limits_L {\left( {{\bf{i}}W_x  + {\bf{j}}W_y  + {\bf{k}}W_z } \right)}  \\
 \cdot \left( {{\bf{i}}x^{\mu  - 1} dx + {\bf{j}}y^{\mu  - 1} dy + {\bf{k}}z^{\mu  - 1} dz} \right) \\
  = \mu \oint\limits_L {W_x x^{\mu  - 1} dx + W_y y^{\mu  - 1} dy + W_z z^{\mu  - 1} dz} . \\
 \end{array}
\end{equation}

From (\ref{eq76}) the Stokes-like theorem for the fractal vector field ${\rm {\bf
W}}$ states that
\begin{equation}
\label{eq82}
\mathop{\iint} \limits_S {\left( {\nabla ^\mu \times {\rm {\bf W}}}
\right)\cdot {\rm {\bf n}}dS} =\oint\limits_L {{\rm {\bf W}}\cdot d{\rm {\bf
l}}} ,
\end{equation}
where ${\rm {\bf W}}$ is a Hausdorff differentiable vector field,
$S$ denotes an open, two sided curve surface, $L$ represents the closed
contour bounding $S$, and $d{\rm {\bf l}}$ denotes the element of the vector
line.

Taking $d{\rm {\bf S}}={\rm {\bf n}}dS$, we have from (\ref{eq52}) that
\begin{equation}
\label{eq83}
\mathop{\iint} \limits_S {\left( {\nabla ^\mu \times {\rm {\bf W}}}
\right)\cdot d{\rm {\bf S}}} =\oint\limits_L {{\rm {\bf W}}\cdot d{\rm {\bf
l}}} .
\end{equation}
It is shown that (\ref{eq83}) becomes the Stokes theorem due to Stokes \cite{19} when $\mu
=1$.

\begin{itemize}
\item \textbf{The Green-like theorem for the fractal vector field}
\end{itemize}
Let us consider
\begin{equation}
\label{eq84}
\begin{array}{l}
 \nabla ^\mu \times {\rm {\bf T}}\\
 =\left( {{\begin{array}{*{20}c}
 {\rm {\bf i}} \hfill & {\rm {\bf j}} \hfill & {\rm {\bf k}} \hfill \\
 {\mu x^{\mu -1}{ }^C\partial _x^{\left( 1 \right)} } \hfill & {\mu y^{\mu
-1}{ }^C\partial _y^{\left( 1 \right)} } \hfill & 0 \hfill \\
 {T_x } \hfill & {T_y } \hfill & \mbox{0} \hfill \\
\end{array} }} \right) \\
 ={\rm {\bf k}}\left( {\mu x^{\mu -1}{ }^C\partial _x^{\left( 1 \right)} T_y
-\mu y^{\mu -1}{ }^C\partial _y^{\left( 1 \right)} T_x } \right) \\
 ={\rm {\bf k}}\mu \left( {x^{\mu -1}{ }^C\partial _x^{\left( 1 \right)} T_y
-y^{\mu -1}{ }^C\partial _y^{\left( 1 \right)} T_x } \right) \\
 \end{array}
\end{equation}
where
\[
{\rm {\bf T}}={\rm {\bf T}}\left( {x,y} \right)=\widetilde{{\rm {\bf
T}}}\left( {x^\mu ,y^\mu } \right)={\rm {\bf i}}T_x +{\rm {\bf j}}T_y .
\]
Hence, by (\ref{eq84}) we show that
\begin{equation}
\label{eq85}
\begin{array}{l}
\mathop{\iint} \limits_S {\left( {\nabla ^\mu \times {\rm {\bf T}}}
\right)\cdot d{\rm {\bf S}}}\\
=
\mu \int\!\!\!\int\limits_S {\left( {x^{\mu -1}{ }^C\partial
_x^{\left( 1 \right)} T_y -y^{\mu -1}{ }^C\partial _y^{\left( 1 \right)} T_x
} \right)dS}, \\
 \end{array}
\end{equation}
where
\[
\begin{array}{l}
\mathop{\iint} \limits_S {\left( {\nabla ^\mu \times {\rm {\bf T}}}
\right)\cdot d{\rm {\bf S}}} \\
=
\mathop{\iint} \limits_S {\left[ {\mu {\rm {\bf
k}}\left( {x^{\mu -1}{ }^C\partial _x^{\left( 1 \right)} T_y -x^{\mu -1}{
}^C\partial _y^{\left( 1 \right)} T_x } \right)} \right]} \cdot \\
 \left[ {{\bf{i}}\mu ^2 y^{\mu  - 1} z^{\mu  - 1} dydz + {\bf{j}}\mu ^2 x^{\mu  - 1} z^{\mu  - 1} dxdz} \right. \\
  + \left. {{\bf{k}}x^{\mu  - 1} y^{\mu  - 1} dxdy} \right]\\
  =
\mu \int\!\!\!\int\limits_S {\left( {x^{\mu -1}{ }^C\partial
_x^{\left( 1 \right)} T_y -y^{\mu -1}{ }^C\partial _y^{\left( 1 \right)} T_x
} \right)dS}, \\
 \end{array}
\]

\[
dS=\mu ^2x^{\mu -1}y^{\mu -1}dxdy,
\]
and
\[
\begin{array}{l}
 d{\bf{S}} = {\bf{i}}\mu ^2 y^{\mu  - 1} z^{\mu  - 1} dydz + {\bf{j}}\mu ^2 x^{\mu  - 1} z^{\mu  - 1} dxdz \\
  + {\bf{k}}\mu ^2 x^{\mu  - 1} y^{\mu  - 1} dxdy. \\
 \end{array}
\]

Moreover, we have
\begin{equation}
\label{eq86}
\begin{array}{l}
 \oint\limits_L {{\rm {\bf T}}\cdot d{\rm {\bf l}}} \\
 =\oint\limits_L {\left(
{{\rm {\bf i}}T_x +{\rm {\bf j}}T_y } \right)\cdot \left[ {\mu \left( {{\rm
{\bf i}}x^{\mu -1}dx+{\rm {\bf j}}y^{\mu -1}dy} \right)} \right]} \\
 =\mu \left( {\oint\limits_L {T_x x^{\mu -1}dx+T_y y^{\mu -1}dy} } \right)
\\
 =\oint\limits_L {T_x dl_x +T_y dl_y }, \\
 \end{array}
\end{equation}
in which
\[
d{\rm {\bf l}}=\mu \left( {{\rm {\bf i}}x^{\mu -1}dx+{\rm {\bf j}}y^{\mu
-1}dy} \right)={\rm {\bf i}}dl_x +{\rm {\bf i}}dl_y ,
\]
where
\[
dl_x =\mu x^{\mu -1}dx
\]
and
\[
dl_y =\mu y^{\mu -1}dy.
\]
From (\ref{eq83}), (\ref{eq85}) and (\ref{eq86}) the Green-like theorem for the fractal vector field
${\rm {\bf {\rm T}}}$ states
\begin{equation}
\label{eq87}
\begin{array}{l}
\mu \left( {\oint\limits_L {T_x x^{\mu -1}dx+T_y y^{\mu -1}dy} } \right)\\
=\mu
^3\int\!\!\!\int\limits_S {\left( {x^{\mu -1}{ }^C\partial _x^{\left( 1
\right)} T_y -y^{\mu -1}{ }^C\partial _y^{\left( 1 \right)} T_x }
\right)x^{\mu -1}y^{\mu -1}dxdy} ,
\end{array}
\end{equation}
which is equal to
\begin{equation}
\label{eq88}
\begin{array}{l}
\oint\limits_L {T_x dl_x +T_y dl_y } \\
=\mu \int\!\!\!\int\limits_S {\left(
{x^{\mu -1}{ }^C\partial _x^{\left( 1 \right)} T_y -y^{\mu -1}{ }^C\partial
_y^{\left( 1 \right)} T_x } \right)dS} ,
\end{array}
\end{equation}
where $S$ is the domain bounded by the contour $L$, and ${\rm {\bf {\rm
T}}}=T_x {\rm {\bf i}}+T_y {\rm {\bf j}}$.

When $\mu =1$, (\ref{eq88}) leads to the Green theorem \cite{16} due to Green \cite{19}.

\begin{itemize}
\item \textbf{The Green-like identities for the fractal vector field}
\end{itemize}
Taking ${\rm {\bf X}}=\Theta \nabla ^\mu \psi$ and ${\rm {\bf Y}}=\psi
\nabla ^\mu \Theta $, we show that
\begin{equation}
\label{eq89}
\begin{array}{l}
\nabla ^\mu \cdot {\rm {\bf X}}=\nabla ^\mu \cdot \left( {\Theta \nabla ^\mu
\psi } \right)=\Theta \nabla ^{2\mu }\psi +\nabla ^\mu \psi \cdot \nabla
^\mu \Theta
\end{array}
\end{equation}
and
\begin{equation}
\label{eq90}
\begin{array}{l}
\nabla ^\mu \cdot {\rm {\bf Y}}=\nabla ^\mu \cdot \left( {\psi \nabla ^\mu
\Theta } \right)=\psi \nabla ^{2\mu }\Theta +\nabla ^\mu \psi \cdot \nabla
^\mu \Theta ,
\end{array}
\end{equation}
where $\psi $ and $\Theta $ are the fractal scalar fields.

With the use of (\ref{eq80}), we may have
\begin{equation}
\label{eq91}
\mathop{\iiint} \limits_{\kern-5.5pt \Omega } {\nabla ^\mu
\cdot {\rm {\bf X}}dV} =\oiint \limits_{\rm {\bf S}} {{\rm {\bf X}}\cdot d{\rm {\bf S}}}
\end{equation}
and
\begin{equation}
\label{eq92}
\mathop{\iint} \limits_{\kern-5.5pt \Omega } {\nabla ^\mu
\cdot {\rm {\bf Y}}dV} =\oiint \limits_{\rm {\bf S}} {{\rm {\bf Y}}\cdot d{\rm {\bf S}}} .
\end{equation}
By using (\ref{eq89}) and (\ref{eq91}) we have that
\begin{equation}
\label{eq93}
\begin{array}{l}
\mathop{\iiint} \limits_{\kern-5.5pt \Omega } {\nabla ^\mu
\cdot \left( {\Theta \nabla ^\mu \psi } \right)dV}\\
=\mathop{\iiint} \limits_{\kern-5.5pt \Omega } {\left(
{\Theta \nabla ^{2\mu }\psi +\nabla ^\mu \psi \cdot \nabla ^\mu \Theta }
\right)dV} \\
=\oiint \limits_{\rm
{\bf S}} {\left( {\Theta \nabla ^\mu \psi } \right)\cdot d{\rm {\bf S}}} .
\end{array}
\end{equation}
By using (\ref{eq90}) and (\ref{eq92}) we show that
\begin{equation}
\label{eq94}
\begin{array}{l}
\mathop{\iiint} \limits_{\kern-5.5pt \Omega } {\nabla ^\mu
\cdot \left( {\psi \nabla ^\mu \Theta } \right)dV}\\
=\mathop{\iiint} \limits_{\kern-5.5pt \Omega } {\left(
{\psi \nabla ^{2\mu }\Theta +\nabla ^\mu \psi \cdot \nabla ^\mu \Theta }
\right)dV} \\
=\oiint \limits_{\rm
{\bf S}} {\left( {\psi \nabla ^\mu \Theta } \right)\cdot d{\rm {\bf S}}} .
\end{array}
\end{equation}
Making use of (\ref{eq93}), we give
\[
\begin{array}{l}
\mathop{\iiint} \limits_{\kern-5.5pt \Omega } {\left(
{\Theta \nabla ^{2\mu }\psi +\nabla ^\mu \psi \cdot \nabla ^\mu \Theta }
\right)dV} =\oiint \limits_{\rm
{\bf S}} {\left( {\Theta \nabla ^\mu \psi } \right)\cdot d{\rm {\bf S}}} ,
\end{array}
\]
which leads to
\begin{equation}
\label{eq95}
\begin{array}{l}
\mathop{\iiint} \limits_{\kern-5.5pt \Omega } {\left(
{\Theta \nabla ^{2\mu }\psi +\nabla ^\mu \psi \cdot \nabla ^\mu \Theta }
\right)dV}
=\oiint \limits_{\rm
{\bf S}} {\Theta \partial _n^\mu \psi dS} ,
\end{array}
\end{equation}
where
\[
\begin{array}{l}
\oiint \limits_{\rm {\bf S}}
{\left( {\Theta \nabla ^\mu \psi } \right)\cdot d{\rm {\bf S}}}\\
=\oiint \limits_{\rm {\bf S}}
{\left( {\Theta \nabla ^\mu \psi } \right)\cdot {\rm {\bf n}}dS}\\
=\oiint \limits_{\rm {\bf S}}
{\Theta \partial _n^\mu \psi dS} .
\end{array}
\]
Here, (\ref{eq95}) is the Green-like identity of first type.

With the aid of (\ref{eq94}) we show
\[
\begin{array}{l}
\mathop{\iiint} \limits_{\kern-5.5pt \Omega } {\left( {\psi
\nabla ^{2\mu }\Theta +\nabla ^\mu \psi \cdot \nabla ^\mu \Theta }
\right)dV}
=\oiint \limits_{\rm
{\bf S}} {\left( {\psi \nabla ^\mu \Theta } \right)\cdot d{\rm {\bf S}}} ,
\end{array}
\]
which leads to
\begin{equation}
\label{eq96}
\begin{array}{l}
\mathop{\iiint} \limits_{\kern-5.5pt \Omega } {\left( {\psi
\nabla ^{2\mu }\Theta +\nabla ^\mu \psi \cdot \nabla ^\mu \Theta }
\right)dV}
=\oiint \limits_{\rm
{\bf S}} {\psi \partial _n^\mu \Theta dS} ,
\end{array}
\end{equation}
where
\[
\begin{array}{l}
\oiint \limits_{\rm {\bf S}}
{\left( {\psi \nabla ^\mu \Theta } \right)\cdot d{\rm {\bf S}}}\\
=\oiint \limits_{\rm {\bf S}}
{\left( {\psi \nabla ^\mu \Theta } \right)\cdot {\rm {\bf n}}dS}\\
=\oiint \limits_{\rm {\bf S}}
{\psi \partial _n^\mu \Theta dS} .
\end{array}
\]
By using (\ref{eq95}) and (\ref{eq96}), we give
\[
\begin{array}{l}
\mathop{\iiint} \limits_{\kern-5.5pt \Omega } {\nabla ^\mu
\cdot \left( {\Theta \nabla ^{2\mu }\psi -\psi \nabla ^{2\mu }\Theta }
\right)dV} \\
=\oiint \limits_{\rm
{\bf S}} {\left( {\Theta \nabla ^\mu \psi -\psi \nabla ^\mu \Theta }
\right)\cdot d{\rm {\bf S}}}
\end{array}
\]
which implies that
\begin{equation}
\label{eq97}
\begin{array}{l}
\mathop{\iiint} \limits_{\kern-5.5pt \Omega } {\nabla ^\mu
\cdot \left( {\Theta \nabla ^{2\mu }\psi -\psi \nabla ^{2\mu }\Theta }
\right)dV} \\
=\oiint \limits_{\rm
{\bf S}} {\left( {\Theta \partial _n^{ \mu } \psi -\psi
\partial _n^{\mu } \Theta } \right)dS} .
\end{array}
\end{equation}
Here, (\ref{eq97}) is the Green-like identity of second type.

When $\mu =1$, Eqs. (\ref{eq95}) and (\ref{eq97}) yield the Green-like identities \cite{16} due
to Green \cite{20}.

The theory of the Hausdorff vector calculus is the special case of the
theory of the general vector calculus \cite{13,14,15} by using the idea of the
Gibbs's expression \cite{21}.

\section{Modelling the fractal power-law flow}\label{Sec:4}

We now consider the coordinate system, expressed in the form:
\[
\left( {t^\mu ,x^\mu ,y^\mu ,z^\mu } \right)=t^\mu +{\rm {\bf i}}x^\mu +{\rm
{\bf j}}y^\mu +{\rm {\bf k}}z^\mu
\]
where ${\rm {\bf i}}$, ${\rm {\bf j}}$ and ${\rm {\bf k}}$ are the unit
vectors in the Cartesian coordinate system.

\begin{itemize}
\item \textbf{The material Hausdorff derivative}
\end{itemize}
Let us consider
\[
\varphi =\varphi \left( {t,x,y,z} \right)=\widetilde{\varphi }\left(
{t,x^\mu ,y^\mu ,z^\mu } \right)
\]
be the fractal fluid field.

The total Hausdorff differential of the fractal scalar field $\varphi $ is
expressed by:
\begin{equation}
\label{eq98}
\begin{array}{l}
 d\varphi
 =\left( {\mu x^{\mu -1}{ }^C\partial _x^{\left( 1 \right)} \varphi
} \right)dx+\left( {\mu y^{\mu -1}{ }^C\partial _y^{\left( 1 \right)}
\varphi } \right)dy\\
+\left( {\mu z^{\mu -1}{ }^C\partial _z^{\left( 1
\right)} \varphi } \right)dz+\partial _t^{\left( 1 \right)} \varphi dt, \\
 \end{array}
\end{equation}
which implies that
\begin{equation}
\label{eq99}
\begin{array}{l}
 \frac{d\varphi }{dt}
 =\left( {\mu x^{\mu -1}{ }^C\partial _x^{\left( 1
\right)} \varphi } \right)\frac{\partial x}{\partial t}+\left( {\mu y^{\mu
-1}{ }^C\partial _y^{\left( 1 \right)} \varphi } \right)\frac{\partial
y}{\partial t}\\
+\left( {\mu z^{\mu -1}{ }^C\partial _z^{\left( 1 \right)}
\varphi } \right)\frac{\partial z}{\partial t}+\partial _t^{\left( 1
\right)} \varphi. \\
 \end{array}
\end{equation}
The material Hausdorff derivative for the fractal fluid
density $\varphi $ is defined as \cite{15}:
\begin{equation}
\label{eq100}
\frac{D\varphi }{Dt}=\partial _t^{\left( 1 \right)} \varphi +{\rm {\bf
\upsilon }}\cdot \nabla ^\mu \varphi
\end{equation}
where
\[
\begin{array}{l}
\nabla ^\mu \varphi
 =\mu \left( {{\rm {\bf i}}x^{\mu -1}{ }^C\partial
_x^{\left( 1 \right)} \varphi +{\rm {\bf j}}y^{\mu -1}{ }^C\partial
_y^{\left( 1 \right)} \varphi +{\rm {\bf k}}z^{\mu -1}{ }^C\partial
_z^{\left( 1 \right)} \varphi } \right)
\end{array}
\]
and ${\rm {\bf \upsilon }}=\left( {\partial x/\partial t,\partial y/\partial
t,\partial z/\partial t} \right)=i\upsilon _x +j\upsilon _y +k\upsilon _z $
is the velocity vector.

When $\mu =1$, Eq. (\ref{eq100}) leads to the material derivative due to Stokes \cite{22}.

\begin{itemize}
\item \textbf{The transport theorem for the fractal power-law fluid}
\end{itemize}
By using (\ref{eq100}), the transport theorem for the fractal power-law fluid $G$ is
expressed as
\begin{equation}
\label{eq101}
\begin{array}{l}
\frac{D}{Dt}\mathop{\iiint} \limits_{\kern-5.5pt {\Omega
\left( t \right)}} {GdV}
=\mathop{\iiint} \limits_{\kern-5.5pt {\Omega \left( t
\right)}} {\left( {\partial _t^{\left( 1 \right)} G+{\rm {\bf \upsilon
}}\cdot \nabla ^\mu G} \right)dV} ,
\end{array}
\end{equation}
which yields that
\begin{equation}
\label{eq102}
\begin{array}{l}
\frac{D}{Dt}\mathop{\iiint} \limits_{\kern-5.5pt {\Omega
\left( t \right)}} {GdV}
=\mathop{\iiint} \limits_{\kern-5.5pt {\Omega \left( t
\right)}} {\partial _t^{\left( 1 \right)} GdV}
+\oiint \limits_{{\rm {\bf
S}}\left( t \right)} {G{\rm {\bf \upsilon }}\cdot d{\rm {\bf S}}}
\end{array}
\end{equation}
because
\begin{equation}
\label{eq103}
\begin{array}{l}
\mathop{\iiint} \limits_{\kern-5.5pt {\Omega \left( t
\right)}} {{\rm {\bf \upsilon }}\cdot \nabla ^\mu GdV}
=\oiint \limits_{{\rm {\bf
S}}\left( t \right)} {G\left( {{\rm {\bf \upsilon }}\cdot {\rm {\bf n}}}
\right)dS}
=\oiint \limits_{{\rm
{\bf S}}\left( t \right)} {G{\rm {\bf \upsilon }}\cdot d{\rm {\bf S}}} ,
\end{array}
\end{equation}
where ${\rm {\bf S}}\left( t \right)$ is the surface of $\Omega \left( t
\right)$, ${\rm {\bf n}}$ is the unit normal to the surface, ${\rm {\bf
\upsilon }}$ is the velocity vector, and $G=G\left( {t,x,y,z}
\right)=\widetilde{G}\left( {t,x^\mu ,y^\mu ,z^\mu } \right)$ is the fractal
power-law fluid.

When $\mu =1$, Eqs. (\ref{eq101}) and (\ref{eq102}) yield the Reynolds transport theorem due
to Reynolds \cite{23}.

\begin{itemize}
\item \textbf{The conservation of the mass for the fractal power-law fluid}
\end{itemize}
The mass of the fractal power-law fluid is defined as
\begin{equation}
\label{eq104}
\begin{array}{l}
\mathop{\iiint} \limits_{\kern-5.5pt {\Omega \left( t
\right)}} {\rho dV} ={\rm M}
\end{array}
\end{equation}
where $\rho =\rho \left( {t,x,y,z} \right)=\widetilde{\rho }\left( {t,x^\mu
,y^\mu ,z^\mu } \right)$ and ${\rm M}={\rm M}\left( {t,x,y,z}
\right)=\widetilde{{\rm M}}\left( {t,x^\mu ,y^\mu ,z^\mu } \right)$.

The conservation of the mass for the fractal power-law fluid is represented
in the form:
\begin{equation}
\label{eq105}
\partial _t^{\left( 1 \right)} \rho +{\rm {\bf \upsilon }}\cdot \nabla ^\mu
\rho =0
\end{equation}
because
\begin{equation}
\label{eq106}
\begin{array}{l}
\frac{D}{Dt} \mathop{\iiint} \limits_{\kern-5.5pt {\Omega
\left( t \right)}} {\rho dV}
=\mathop{\iiint} \limits_{\kern-5.5pt {\Omega \left( t
\right)}} {\left( {\partial _t^{\left( 1 \right)} \rho +{\rm {\bf \upsilon
}}\cdot \nabla ^\mu \rho } \right)dV} =0,
\end{array}
\end{equation}
where ${\rm {\bf \upsilon }}$ is the velocity vector.

Let ${\rm {\bf \upsilon }}$ be a constant. Then we have
\begin{equation}
\label{eq107}
\partial _t^{\left( 1 \right)} \rho +\nabla ^\mu \cdot \left( {{\rm {\bf
\upsilon }}\rho } \right)=0
\end{equation}
which is derived from
\begin{equation}
\label{eq108}
\begin{array}{l}
\frac{D}{Dt}\mathop{\iiint} \limits_{\kern-5.5pt {\Omega
\left( t \right)}} {\rho dV}
=\mathop{\iiint} \limits_{\kern-5.5pt {\Omega \left( t
\right)}} {\left[ {\partial _t^{\left( 1 \right)} \rho +\nabla ^\mu \left(
{{\rm {\bf \upsilon }}\cdot \rho } \right)} \right]dV} =0.
\end{array}
\end{equation}
When $\mu =1$, Eqs. (\ref{eq105}) and (\ref{eq106}) are the expressions for the classical
conservation of the mass due to Euler \cite{24} and Lagrange \cite{25}.

\begin{itemize}
\item \textbf{The velocity gradient tensor for the fractal power-law fluid}
\end{itemize}
We now consider the velocity gradient tensor for the fractal power-law
fluid, expressed in the form
\begin{equation}
\label{eq109}
\nabla ^\mu \cdot {\rm {\bf \upsilon }}=\frac{1}{2}\left( {\varsigma +\tau }
\right)+\frac{1}{2}\left( {\varsigma -\tau } \right)={\rm {\bf \eta
}}+\frac{1}{2}\left( {\varsigma -\tau } \right)
\end{equation}
in which
\begin{equation}
\label{eq110}
\varsigma =\nabla ^{\left( {D_1 ,D_2 ,D_3 } \right)}\cdot {\rm {\bf \upsilon
}}=0,
\end{equation}
where the strain tensor for the fractal power-law fluid is defined as
\[
{\rm {\bf \eta }}=\left( {\varsigma +\tau } \right)/2
\]
with velocity gradient $\varsigma =\nabla ^\mu \cdot {\rm {\bf \upsilon }}$
and $\tau ={\rm {\bf \upsilon }}\cdot \nabla ^\mu $.

The stress tensor for the fractal power-law fluid is defined as
\begin{equation}
\label{eq111}
{\rm {\bf {\rm H}}}=-p{\rm {\bf I}}+2\varepsilon {\rm {\bf \eta }},
\end{equation}
where $\varepsilon $ are the shear moduli of viscosity, and ${\rm {\bf
I}}$is the unit tensor.

When$\mu =1$, the velocity gradient tensor for the fractal power-law fluid
implies the Cauchy strain tensor due to Cauchy \cite{26}.

\begin{itemize}
\item \textbf{The conservation of the momentums for the fractal power-law fluid}
\end{itemize}
The conservation of the momentums for the fractal power-law fluid reads
\begin{equation}
\label{eq112}
\begin{array}{l}
\frac{D}{Dt}\mathop{\iiint} \limits_{\kern-5.5pt {\Omega
\left( t \right)}} {\rho {\rm {\bf \upsilon }}dV}
=\mathop{\iiint} \limits_{\kern-5.5pt {\Omega \left( t
\right)}} {{\rm {\bf b}}dV} +\oiint \limits_{{\rm {\bf S}}\left( t \right)} {{\rm {\bf {\rm H}}}\cdot
d{\rm {\bf S}}} ,
\end{array}
\end{equation}
where ${\rm {\bf b}}$ is the specific body force.

By using (\ref{eq101}) and (\ref{eq102}), we have
\begin{equation}
\label{eq113}
\begin{array}{l}
\frac{D}{Dt}\mathop{\iiint} \limits_{\kern-5.5pt {\Omega
\left( t \right)}} {\rho {\rm {\bf \upsilon }}dV}
=\mathop{\iiint} \limits_{\kern-5.5pt {\Omega \left( t
\right)}} {\partial _t^{\left( 1 \right)} \left( {\rho {\rm {\bf \upsilon
}}} \right)dV} +\oiint \limits_{{\rm {\bf S}}\left( t \right)} {\left( {\rho {\rm {\bf
\upsilon }}} \right){\rm {\bf \upsilon }}\cdot d{\rm {\bf S}}}
\end{array}
\end{equation}
and
\begin{equation}
\label{eq114}
\begin{array}{l}
\frac{D}{Dt}\mathop{\iiint} \limits_{\kern-5.5pt {\Omega
\left( t \right)}} {\rho {\rm {\bf \upsilon }}dV}
=\mathop{\iiint} \limits_{\kern-5.5pt {\Omega \left( t
\right)}} {\left[ {\partial _t^{\left( 1 \right)} \left( {\rho {\rm {\bf
\upsilon }}} \right)+{\rm {\bf \upsilon }}\cdot \nabla ^\mu \left( {\rho
{\rm {\bf \upsilon }}} \right)} \right]dV} .
\end{array}
\end{equation}
With (\ref{eq80}) we present
\begin{equation}
\label{eq115}
\begin{array}{l}
\oiint \limits_{{\rm {\bf
S}}\left( t \right)} {{\rm {\bf {\rm H}}}\cdot d{\rm {\bf S}}}
=\mathop{\iiint} \limits_{\kern-5.5pt {\Omega \left( t
\right)}} {\nabla ^\mu \cdot {\rm {\bf {\rm H}}}dV} .
\end{array}
\end{equation}
From (\ref{eq112}) we have
\[
\begin{array}{l}
\mathop{\iiint} \limits_{\kern-5.5pt {\Omega \left( t
\right)}} {\left[ {\partial _t^{\left( 1 \right)} \left( {\rho {\rm {\bf
\upsilon }}} \right)+{\rm {\bf \upsilon }}\cdot \nabla ^\mu \left( {\rho
{\rm {\bf \upsilon }}} \right)} \right]dV}
=\mathop{\iiint} \limits_{\kern-5.5pt {\Omega \left( t
\right)}} {{\rm {\bf b}}dV}
+\mathop{\iiint} \limits_{\kern-5.5pt {\Omega \left( t
\right)}} {\nabla ^\mu \cdot {\rm {\bf {\rm H}}}dV} ,
\end{array}
\]
which leads to
\begin{equation}
\label{eq116}
\begin{array}{l}
\mathop{\iiint} \limits_{\kern-5.5pt {\Omega \left( t
\right)}} {\left[ {\partial _t^{\left( 1 \right)} \left( {\rho {\rm {\bf
\upsilon }}} \right)+{\rm {\bf \upsilon }}\cdot \nabla ^\mu \left( {\rho
{\rm {\bf \upsilon }}} \right)-\nabla ^\mu \cdot {\rm {\bf {\rm H}}}-{\rm
{\bf b}}} \right]dV}
=0.
\end{array}
\end{equation}
From (\ref{eq116}) it follows that
\begin{equation}
\label{eq117}
\begin{array}{l}
\partial _t^{\left( 1 \right)} \left( {\rho {\rm {\bf \upsilon }}}
\right)+{\rm {\bf \upsilon }}\cdot \nabla ^\mu \left( {\rho {\rm {\bf
\upsilon }}} \right)-\nabla ^\mu \cdot {\rm {\bf {\rm H}}}-{\rm {\bf b}}=0,
\end{array}
\end{equation}
which leads to
\begin{equation}
\label{eq118}
\begin{array}{l}
\rho \left( {\partial _t^{\left( 1 \right)} {\rm {\bf \upsilon }}+{\rm {\bf
\upsilon }}\cdot \nabla ^\mu {\rm {\bf \upsilon }}} \right)-\nabla ^\mu
\cdot {\rm {\bf {\rm H}}}-{\rm {\bf b}}=0.
\end{array}
\end{equation}
From (\ref{eq111}) and (\ref{eq118}) we show that
\begin{equation}
\label{eq119}
\begin{array}{l}
\nabla ^\mu \cdot {\rm {\bf {\rm H}}}=-\nabla ^\mu p+\varepsilon \nabla
^{2\mu }{\rm {\bf \upsilon }}
\end{array}
\end{equation}
such that
\begin{equation}
\label{eq120}
\begin{array}{l}
\rho \left( {\partial _t^{\left( 1 \right)} {\rm {\bf \upsilon }}+{\rm {\bf
\upsilon }}\cdot \nabla ^\mu {\rm {\bf \upsilon }}} \right)=-\nabla ^\mu
p+\varepsilon \nabla ^{2\mu }{\rm {\bf \upsilon }}+{\rm {\bf b}}.
\end{array}
\end{equation}
From (\ref{eq120}) we get the system of the power-law flow, given by the
the fractal power-law equations
\begin{equation}
\label{eq121}
\begin{array}{l}
\rho \left( {\partial _t^{\left( 1 \right)} {\rm {\bf \upsilon }}+{\rm {\bf
\upsilon }}\cdot \nabla ^\mu {\rm {\bf \upsilon }}} \right)=-\nabla ^\mu
p+\varepsilon \nabla ^{2\mu }{\rm {\bf \upsilon }}+{\rm {\bf b}},
\end{array}
\end{equation}
and
\begin{equation}
\label{eq122}
\begin{array}{l}
\nabla ^\mu \cdot {\rm {\bf \upsilon }}=0.
\end{array}
\end{equation}
When $\mu =1$, Eqs. (\ref{eq121}) and (\ref{eq122}) become the system of the Navier-Stokes
equations for the fluid due to Navier \cite{27} and Stokes \cite{22}. The power-law
flow in the real world problems is the special case of the theory of the
general and power-law flow in the real world problems studied in \cite{13,14,15}.

Let ${\rm F}={\rm F}\left( {t,x,y,z}
\right)=\widetilde{{\rm F}}\left( {t,x^\mu ,y^\mu ,z^\mu } \right)$ be a curve.
If there exist all Chen Hausdorff derivatives of the curve ${\rm F}$ in the space domain ${{x^\mu}\times{y^\mu}\times{z^\mu}}\in {\mathbb{R}^{3\mu}}$
and all derivatives of the curve ${\rm F}$ in the time domain $t \in \mathbb{R}$,
then ${\rm F}$ is smooth in the space-time domain, where
$0<\mu \leq 1$.

\textbf{Conjecture}
Do the fractal power-law equations (\ref{eq121}) and (\ref{eq122}) on a fractal domain
${\Omega \left( t\right)}$ in ${\mathbb{R}^{3\mu}}$ have a unique smooth solution for all time $t\geq0$?

It is easy to see that \textbf{Conjecture} is analogous to the Smale¨s 15th Problem \cite{28}.
When we take $\mu =1$, \textbf{Conjecture} becomes the Smale¨s 15th Problem \cite{28} or one of the
Millennium Prize Problems for the Navier--Stokes equations \cite{29}.

\section{Modelling the anomalous diffusion equation} \label{Sec:5}

When $p$ is a constant, we have from (\ref{eq119}) that
\begin{equation}
\label{eq123}
\nabla ^\mu \cdot {\rm {\bf {\rm H}}}=\varepsilon \nabla ^{2\mu }{\rm {\bf
\upsilon }},
\end{equation}
which implies that
\begin{equation}
\label{eq124}
\begin{array}{l}
\oiint \limits_{{\rm {\bf
S}}\left( t \right)} {{\rm {\bf {\rm H}}}\cdot d{\rm {\bf S}}}\\
=\mathop{\iiint} \limits_{\kern-5.5pt {\Omega \left( t
\right)}} {\nabla ^\mu \cdot {\rm {\bf {\rm H}}}dV} \\
=\varepsilon
\mathop{\int\!\!\!\int\!\!\!\int}\limits_{\kern-5.5pt {\Omega \left( t
\right)}} {\nabla ^{2\mu }{\rm {\bf \upsilon }}dV}\\
=\varepsilon
\mathop{\iiint} \limits_{\kern-5.5pt {\Omega \left( t
\right)}} {\nabla ^\mu \cdot \left( {\nabla ^\mu {\rm {\bf \upsilon }}}
\right)dV} \\
=\varepsilon \oiint \limits_{\rm {\bf S}} {\left( {\nabla ^\mu {\rm {\bf \upsilon }}}
\right)\cdot d{\rm {\bf S}}} .
\end{array}
\end{equation}
Thus, from (\ref{eq124}), we get
\begin{equation}
\label{eq125}
\begin{array}{l}
\oiint \limits_{{\rm {\bf
S}}\left( t \right)} {{\rm {\bf {\rm H}}}\cdot d{\rm {\bf S}}}\\
=\mathop{\iiint} \limits_{\kern-5.5pt {\Omega \left( t
\right)}} {\nabla ^\mu \cdot {\rm {\bf {\rm H}}}dV} \\
=\varepsilon
\mathop{\iiint} \limits_{\kern-5.5pt {\Omega \left( t
\right)}} {\nabla ^{2\mu }{\rm {\bf \upsilon }}dV}\\
=\varepsilon
\oiint \limits_{\rm {\bf S}}
{\left( {\nabla ^\mu {\rm {\bf \upsilon }}} \right)\cdot d{\rm {\bf S}}} .
\end{array}
\end{equation}
By using (\ref{eq112}) we show
\begin{equation}
\label{eq126}
\begin{array}{l}
\mathop{\iiint} \limits_{\kern-5.5pt {\Omega \left( t
\right)}} {\left[ {\partial _t^{\left( 1 \right)} \left( {\rho {\rm {\bf
\upsilon }}} \right)+{\rm {\bf \upsilon }}\cdot \nabla ^\mu \left( {\rho
{\rm {\bf \upsilon }}} \right)} \right]dV}\\
=\mathop{\iiint} \limits_{\kern-5.5pt {\Omega \left( t
\right)}} {{\rm {\bf b}}dV}
+\mathop{\iiint} \limits_{\kern-5.5pt {\Omega \left( t
\right)}} {\nabla ^\mu \cdot {\rm {\bf {\rm H}}}dV} \\
 =\mathop{\iiint} \limits_{\kern-5.5pt {\Omega \left( t
\right)}} {{\rm {\bf b}}dV} +\varepsilon
\mathop{\iiint} \limits_{\kern-5.5pt {\Omega \left( t
\right)}} {\nabla ^{2\mu }{\rm {\bf \upsilon }}dV} \\
 =\mathop{\iiint} \limits_{\kern-5.5pt {\Omega \left( t
\right)}} {{\rm {\bf b}}dV} +\varepsilon
\oiint \limits_{\rm {\bf S}}
{\left( {\nabla ^\mu {\rm {\bf \upsilon }}} \right)\cdot d{\rm {\bf S}}} \\
 \end{array}
\end{equation}
which is equal to
\begin{equation}
\label{eq127}
\begin{array}{l}
\mathop{\iiint} \limits_{\kern-5.5pt {\Omega \left( t
\right)}} {\left[ {\partial _t^{\left( 1 \right)} \left( {\rho {\rm {\bf
\upsilon }}} \right)+{\rm {\bf \upsilon }}\cdot \nabla ^\mu \left( {\rho
{\rm {\bf \upsilon }}} \right)} \right]dV}
=\mathop{\iiint} \limits_{\kern-5.5pt {\Omega \left( t
\right)}} {{\rm {\bf b}}dV} +\varepsilon
\mathop{\iiint} \limits_{\kern-5.5pt {\Omega \left( t
\right)}} {\nabla ^{2\mu }{\rm {\bf \upsilon }}dV} .
\end{array}
\end{equation}
From (\ref{eq127}) we have
\begin{equation}
\label{eq128}
\begin{array}{l}
\mathop{\iiint} \limits_{\kern-5.5pt {\Omega \left( t
\right)}} {\left[ {\partial _t^{\left( 1 \right)} \left( {\rho {\rm {\bf
\upsilon }}} \right)+{\rm {\bf \upsilon }}\cdot \nabla ^\mu \left( {\rho
{\rm {\bf \upsilon }}} \right)-{\rm {\bf b}}-\varepsilon \nabla ^{2\mu }{\rm
{\bf \upsilon }}} \right]dV}\\
 =0
\end{array}
\end{equation}
such that
\begin{equation}
\label{eq129}
\partial _t^{\left( 1 \right)} \left( {\rho {\rm {\bf \upsilon }}}
\right)+{\rm {\bf \upsilon }}\cdot \nabla ^\mu \left( {\rho {\rm {\bf
\upsilon }}} \right)-\varepsilon \nabla ^{2\mu }{\rm {\bf \upsilon }}={\rm
{\bf b}}.
\end{equation}
When ${\rm {\bf b}}={\rm {\bf 0}}$, we show that
\begin{equation}
\label{eq130}
\partial _t^{\left( 1 \right)} \left( {\rho {\rm {\bf \upsilon }}}
\right)+{\rm {\bf \upsilon }}\cdot \nabla ^\mu \left( {\rho {\rm {\bf
\upsilon }}} \right)=\varepsilon \nabla ^{2\mu }{\rm {\bf \upsilon }},
\end{equation}
which implies that
\begin{equation}
\label{eq131}
\partial _t^{\left( 1 \right)} {\rm {\bf \upsilon }}+{\rm {\bf \upsilon
}}\cdot \nabla ^\mu {\rm {\bf \upsilon }}=\vartheta \nabla ^{2\mu }{\rm {\bf
\upsilon }},
\end{equation}
where $\vartheta =\varepsilon /\rho $.

In one-dimensional case, we have from (\ref{eq131}) that
\begin{equation}
\label{eq132}
\partial _t^{\left( 1 \right)} \upsilon +\upsilon \mu x^{\mu -1}{
}^C\partial _x^{\left( 1 \right)} \upsilon =\vartheta \mu ^2x^{2\mu -2}{
}^C\partial _x^{\left( 2 \right)} \upsilon ,
\end{equation}
where $\upsilon =\upsilon \left( {t,x} \right)=\widetilde{\upsilon }\left(
{t,x^\mu } \right)$.

When we neglect the nonlinear term in (\ref{eq131}), the linear anomalous diffusion
equation is expressed in the form:
\begin{equation}
\label{eq133}
\partial _t^{\left( 1 \right)} \upsilon =\vartheta \mu ^2x^{2\mu -2}{
}^C\partial _x^{\left( 2 \right)} \upsilon ,
\end{equation}
where $\upsilon =\upsilon \left( {t,x} \right)=\widetilde{\upsilon }\left(
{t,x^\mu } \right)$.

Eq. (\ref{eq132}) is an anomalous diffusion equation in the real theory of the
turbulent fluid motion.

When we take $\mu =1$, Eq. (\ref{eq132}) becomes the Burgers diffusion equation due
to Burgers \cite{30}.

\section{Conclusion} \label{Sec:6}

In our work we showed the theory of the Hausdorff vector calculus with use
of the Chen Hausdorff calculus. We discussed the Gauss-Ostrogradsky-like,
Stokes-like, and Green-like theorems, and Green-like identities for the
fractal field. By using the Hausdorff vector calculus we obtained the theory
of the fractal power-law flow analogous to the Navier-Stokes and Burgers
diffusion equations. A conjecture for the fractal power-law flow equations analogous
to the Smale¨s 15th Problem (one of the Millennium Prize Problems
for the Navier--Stokes equations) has been also addressed.
This plays an important role in the study of the
anomalous and complex behaviors of the flows in the real world problems.

\section*{ACKNOWLEDGMENTS}

This work is supported by the Yue-Qi Scholar of the China University of Mining and Technology (No. 102504180004).

\bibliographystyle{fractals}


\end{document}